\documentclass[10pt,twoside]{amsart}

\usepackage{mathtools,amssymb,csquotes}

\newcounter{pnum}[section]
\renewcommand{\thepnum}{\thesection.\arabic{pnum}}
\newcommand{\nump}{\noindent\refstepcounter{pnum}{\scriptsize\thepnum}\quad}

\DeclareMathOperator{\Ann}{Ann}

\newcommand{\field}[1]{\mathbb{#1}}
\newcommand{\C}{\field C}

\newcommand{\R}{\field R}
\newcommand{\Z}{\field Z}

\newcommand{\Lie}[1]{\mathit{#1}}
\newcommand{\lie}[1]{\mathfrak{#1}}

\newcommand{\SU}[1]{\Lie SU(#1)}

\newcommand{\hook}{\mathop{\lrcorner}}

\DeclarePairedDelimiterX{\Set}[1]\{\}{%
  
#1 }
\DeclarePairedDelimiterX{\Span}[1]{\langle}{\rangle}{%
  
#1 }

\begin{document}

\title{Solvable groups and a shear construction}

\author{Marco Freibert}
\address{Department of Mathematics\\
Aarhus University\\
Ny Munkegade 118, Bldg 1530\\
DK-8000 Aarhus C\\
Denmark}
\email{marco.freibert@math.au.dk}
\thanks{Marco Freibert partially supported by Danish
Council for Independent Research \textbar\ Natural Sciences project
DFF - 4002-00125.}
\author{Andrew Swann}
\address{Department of Mathematics\\
Aarhus University\\
Ny Munkegade 118, Bldg 1530\\
DK-8000 Aarhus C\\
Denmark}
\email{swann@math.au.dk}

\dedicatory{Dedicated to Prof.\ Jaime Mu\~noz-Masqu\'e on the occasion of his
65th birthday}

\begin{abstract}
  \noindent The twist construction is a geometric model of T-duality
  that includes constructions of nilmanifolds from tori.  This paper
  shows how one-dimensional foliations on manifolds may be used in a
  shear construction, which in algebraic form builds certain solvable
  Lie groups from Abelian ones.  We discuss other examples of
  geometric structures that may be obtained from the shear
  construction.
\end{abstract}

\subjclass[2010]{Primary 53C15; Secondary 53C12, 53C30, 53C55, 32C37}

\keywords{T-duality, solvable Lie group, foliation, Hermitian geometry,
G2 structure}

\maketitle

\section{Introduction}
\label{sec:introduction}

Recent years have seen a large number of constructions and
classifications of geometric structures on nilpotent and solvable Lie
groups, for example \cite{Fino-PS:SKT,Ugarte:balanced,%
Conti-FS:qc-closed,Barberis-DV:canonical}. 
Many of these structures are motivated by ideas from theoretical
physics, and particularly various requirements coming from string and
M-theories.  
An important aspect of such theories are various duality relations.
In particular, as Strominger, Yau \& Zaslow
\cite{Strominger-YZ:mirror} proposed that T-duality is closely related
to the concept of mirror Calabi-Yau manifolds.
When fluxes are introduced, the relevant geometries no longer have
special holonomy, but are special types of almost Hermitian manifolds
in dimension~\( 6 \) and, for M-theory, \( G_2 \)-manifolds in
dimension~\( 7 \).

In \cite{Swann:twist,Swann:T} a geometric version of T-duality,
called the twist construction,
was described that reproduces the known results on nilmanifolds and
provides other geometric examples.
It was applied in \cite{Macia-S:c-map} to describe the geometry of the
c-map \cite{Cecotti-FG:II} that constructs quaternionic K\"ahler
manifolds in dimension \( 4n+4 \) from projective special K\"ahler
manifolds in dimension~\( 2n \).
The homogeneous models of this construction
\cite{DeWit-vP:potentials,DeWit-vP:special} provide
all known examples of homogeneous quaternionic K\"ahler metrics on
completely solvable Lie groups \cite{Cortes:Alekseevskian}.
However, the construction of \cite{Macia-S:c-map} requires 
modifying the geometry via so-called elementary deformations,
before the twist construction is used.
This is related to the solvable, rather than nilpotent, nature of the
homogeneous examples.
It is therefore interesting to look for a broader construction that
includes the geometry of solvable groups.  
The purpose of this paper is to propose such a shear construction in
the situation where there is a single symmetry,
or more generally an appropriate one-dimensional foliation.
We will restrict ourselves to a situation which for solvable groups
corresponds to having only real eigenvalues.
In future work, we will describe how foliations from bundles with flat
connections can be used to remove this restriction and consider
foliations of higher rank.

In section~\ref{sec:models}, we consider the underlying Lie algebraic
picture, first recalling the description for nilpotent Lie groups and
then extending it to the solvable case.  
We then provide a general geometric set-up in
section~\ref{sec:geom-constr} based on a double fibration picture, and
discuss carefully what type of bundles may occur.
The fibrations will be seen to be given by principal
bundles with one-dimensional fibres,  
but the connection-like one-forms are not necessarily principal.
In section~\ref{sec:geometric-structures}, we describe how geometric
structures may be transferred through the shear construction,
considering which differential forms are naturally related via the
horizontal distribution,
and providing a formula for the de Rham differentials.
Finally section~\ref{sec:examples} shows how this shear construction
may be applied in examples for different geometric structures.

\section{Algebraic models}
\label{sec:models}

\nump
The model behind the twist construction is based on the geometry of
nilmanifolds.
Recall that a Lie algebra \( \lie n \)
is nilpotent if its lower central series
\( \lie n^{(1)} = \lie n' = [\lie n,\lie n] \),
\( \lie n^{(k)} = [\lie n,\lie n^{(k-1)}] \),
terminates, so \( \lie n^{(r)} = \{0\} \)
for some \( r\geqslant 1 \).
The smallest such number \( r \)
is called the \enquote{step length} of~\( \lie n \).
Dually this condition is that there is a minimal filtration
\( \lie n^* = V_0 > V_1 > \dots > V_{r-1} > \{0\} \),
given by \( V_i = V_i(\lie n) = \Ann(\lie n^{(r-i)}) \),
with \( dV_i \subset \Lambda^2V_{i+1} \)
for each~\( i \).
Our convention is that \( d\alpha(X,Y) = -\alpha([X,Y]) \);
the Jacobi identity is equivalent to \( d\circ d = 0 \).

Fix an element \( \alpha \in V_0 \setminus V_1 \)
and choose a splitting \( \lie n^* = \R\alpha \oplus W \)
with \( V_1 \leqslant W \).
Let \( F \in \Lambda^2V_1 \) be a two-form
with \( dF = 0 \).
Then \( dW \leqslant \Lambda^2V_1 \leqslant \Lambda^2W \)
and we may define a new Lie algebra \( \lie m \),
by taking \( \lie m^* = \R\beta + W \),
retaining the definition of \( d \) on~\( W \)
and putting \( d\beta = d\alpha + F \).
The algebras \( \lie n \) and \( \lie m \) are then said to be
related by a twist construction.

One valid choice in this construction is \( F = -d\alpha \).
This gives \( d\beta = 0 \), and in this case \( \lie m \)
is a Lie algebra direct sum \( \R \oplus W^* \),
with \( W^* \)~nilpotent.
On the other hand, if \( \lie m \) is the twist of \( \lie n \)
via~\( F \), then we may invert the construction by using the
\( 2 \)-form \( -F \).
It follows that any nilpotent algebra may be obtained from repeated
twists of an Abelian algebra of the same dimension.

Note that dual to the splitting \( \lie n^* = \R\alpha \oplus W \)
we get a unique \( X \in \lie n\) specified by \( W(X) = 0 \)
and \( \alpha(X) = 1 \).
The condition \( V_1 \subset W \) ensures that \( X \) is central.

A simple example is provided by the Heisenberg algebra
\( \lie h_3 = (0,0,12) \),
where the abbreviated notation means that \( \lie h^* \) has a basis
\( e_1,e_2,e_3 \) with corresponding differentials \( 0,0,12 \),
meaning \( de_1 = 0 \), \( de_2 = 0 \) and
\( de_3 = e_1\wedge e_2 = e_{12} \).
In this case, we may take \( \alpha = e_3 \) and
\( F = - d\alpha = -e_{12} \).
The resulting twist is the Abelian Lie algebra \( (0,0,0) \).

\medbreak
\nump
Now consider a solvable algebra~\( \lie s \).
This means that the derived series \( \lie s', (\lie s')', \dots \)
terminates at some finite step.
One then has that \( \lie n = \lie s' \) is nilpotent and that
\( \lie a = \lie s/\lie n \) acts on \( \lie n \) as an Abelian
algebra of endomorphisms.
In particular, if \( \lie n \) has step length~\( r \), 
then the subspace \( \lie n^{(r-1)} \) is preserved by~\( \lie a \).
It follows that the complexification \( \lie n^{(r-1)} \otimes \C \)
contains a one-dimensional invariant subspace~\( \xi_\C \).

For the purposes of this article, 
let us work in the case when \( \xi_\C \) may be chosen as
the complexification of a real one-dimensional subspace
\( \xi \leqslant \lie n^{(r-1)} \) preserved by~\( \lie a \).
For example, this will be the case if \( \lie s \) is completely
solvable.

Fix a basis element \( X \) of~\( \xi \).
Let \( \alpha \) be any element of \( \lie s^* \) with
\( \alpha(X) = 1 \), 
then for \( W = \Ann(\xi) \) we have
\( \lie s^* = \R\alpha \oplus W \).

The choice of \( \xi \) implies that there is an element
\( \eta \in \lie s^* \) defined by \( [A,X] = \eta(A)X \) for each
\( A \in \lie s \).
We now have \( d\alpha = \eta \wedge \alpha + F \),
with \( F \in \Lambda^2W \) and \( \eta \in W \).
Furthermore \( \xi \) is an ideal of \( \lie s \)
and \( dW \subset \Lambda^2W \).
Note that the relation \( d^2\alpha = 0 \), implies that \( \eta \)
and \( F \) satisfy
\begin{equation}
  \label{eq:F-eta-s}
  dF = \eta\wedge F,\quad
  d\eta = 0 \quad\text{and}\quad
  \eta|_\xi = 0.
\end{equation}

Now suppose \( F_0 \in \Lambda^2\lie s^* \) is another two-form.
We wish to define a new Lie algebra \( \lie r \) by putting
\( \lie r^* = \R\beta \oplus W \) with
\( d\beta \) to be essentially \( d\alpha + F_0 \).
More precisely, we may write \( F_0 = \eta' \wedge \alpha + F' \),
with \( \eta', F' \in \Lambda^*W \).
Then define
\begin{equation}
  \label{eq:d-beta}
  d\beta = \tilde\eta \wedge \beta + \tilde F,
\end{equation}
with \( \tilde\eta = \eta+\eta' \) and \( \tilde F = F + F' \).
For this to define a Lie algebra, we need \( d^2\beta = 0 \),
which as above is equivalent to
\(d\tilde F = \tilde\eta \wedge \tilde F \)
and \( d\tilde\eta = 0 \).
Translating this back to conditions on~\( F_0 \),
we first note that \( \eta' \) is closed and
\( \eta' = - X \hook F_0 \). 
Now 
\begin{equation*}
  \begin{split}
    dF_0
    &= - \eta'\wedge d\alpha + dF' \\
    &= - \eta'\wedge\eta\wedge\alpha - \eta'\wedge F
      + (\eta + \eta')\wedge(F+F') - \eta\wedge F \\
    &= \eta\wedge\eta'\wedge\alpha + (\eta+\eta')\wedge F'
    = (\eta - X\hook F_0) \wedge F_0.
  \end{split}
\end{equation*}
In other words, we get a Lie algebra if and only if
\begin{equation}
  dF_0 = \eta_0 \wedge F_0,\quad
  d\eta_0 = 0,\quad
  \eta_0|_\xi = 0
\end{equation}
for \( \eta_0 = \eta - X\hook F_0 \).
This is one model of a pair of Lie algebras related by a shear.

A valid choice for \( F_0 \) is \( -F \), which means that
in~\( \lie r \) we have killed some of the Lie brackets in
\( \lie n \) and \( d\beta \) is decomposable.
Another good choice is \( F_0 = -d\alpha \),
which produces a product algebra
\( \lie r = (\lie g/\xi) \oplus \R \).
Conversely, appropriate choices of \( F_0 \)
allows us to construct \( \lie s \) as a shear of \( \lie r \).

In contrast to the twist case above, \( \xi = \Span X \)~is an ideal
which is not necessarily central.

As an example, consider the completely solvable algebra
\( \lie s = (51,52,53,2.54,0) \) which has
Abelian nilradical.
Taking \( \alpha = e_4 \), we have \( F = 0 \) and
\( \eta = 2e_5 \).
The choice \( F_0 = e_{13} \), has
\( dF_0 = e_{513} - e_{153} = \eta\wedge F_0 \),
so \( \eta_0 = \eta \),
and the resulting shear is \( (51,52,53,13+2.54,0) \)
which has non-Abelian nilradical.
Taking \( F_0 = - d\alpha = -2e_{54} \), gives a shear that is
\( \lie r = (51,52,53,0,0) \).
The original algebra \( \lie s \) is obtained from \( \lie r \),
simply by taking \( \alpha = e_4 \) and \( F_0 = 2e_{54} \).

\medbreak
\nump
In both the above constructions, a geometric structure on the original
algebra described by left-invariant tensors may be transferred to a
corresponding geometric structure on the new algebra by replacing each
occurrence of \( \alpha \) by \( \beta \), or vice versa.
It is clear that this will not necessarily preserve integrability
properties of these structures, since these are determined by the
exterior differential.  However, it is straightforward to trace how
such differentials change.

\section{Geometric constructions}
\label{sec:geom-constr}

Let us now turn to geometric versions of the above algebraic
relations.  
A tradition in T-duality has been to consider a space and its dual as
being fibred over a common base.  
In \cite{Swann:twist,Swann:T} a different approach was proposed in
which both spaces are to be regarded as the bases of fibrations from a
common total space.

\smallbreak
\nump
In more detail, let \( M \) be a manifold with a circle action
generated by a vector field~\( X \).
Given a closed two-form with integral periods
\( F \in \Omega^2(M)_\Z \),
there is a principal circle bundle \( P \to M \) and a connection
one-form \( \theta \in \Omega^1(P) \) whose curvature is~\( F \).
If there is a function \( a \in C^\infty(M) \) with
\( da = - X \hook F \), then we may call \( a \) a Hamiltonian,
and define a vector field \( X' \) on~\( P \) preserving \( \theta \) via
\begin{equation*}
  X'  = \widetilde X + a Y.
\end{equation*}
Here \( \widetilde X \in \mathcal H = \ker \theta \) is the horizontal
lift of \( X \), and \( Y \) is the generator of the principal action
on \( P \).
The twist \( W \) of \( M \) by the data \( (X,F,a) \) is then defined
to be
\begin{equation*}
  W = P/\Span{X'},
\end{equation*}
whenever this is a smooth manifold.

For a shear construction we propose the following.  Let \( M \) be a
manifold with a vector field~\( X \).
Write \( \xi \) for the (possibly singular)
distribution~\( \R X \subset TM \).
Suppose \( F \) is a two-form on~\( M \) with the property that
\begin{equation}
  \label{eq:F-eta}
  dF = \eta \wedge F,
\end{equation}
for some closed one-form \( \eta \in \Omega^1(M) \) with
\( \eta|_\xi = 0  \).

Choose a fibre bundle \( P \xrightarrow\pi M \) with connected
one-dimensional fibres and write \( \zeta \) for~\( \ker \pi_* \).  
Suppose \( P \) carries a one-form \( \theta \in \Omega^1(P) \)
with \( \mathcal H = \ker \theta \) transverse to \( \zeta \)
and with
\begin{equation}
  \label{eq:d-theta}
  d\theta = \pi^*\eta \wedge \theta + \pi^*F.
\end{equation}

We now seek a one-dimensional foliation \( \xi' \) of~\( P \)
transverse to \( \mathcal H \) 
with \( \pi_*\xi' = \xi \) and admitting local sections
preserving~\( \theta \).
The shear \( S \) of \( M \) will then be defined to be
the leaf space
\begin{equation*}
  S = P/\Span{\xi'}
\end{equation*}
whenever this is a smooth manifold.

\bigbreak
Let us now discuss various stages in this construction, in order to
justify some of the choices made and to clarify the set-up.

\medbreak
\nump
If the rank of the two-form \( F \) is at least~\( 6 \), then
equation~\eqref{eq:F-eta} implies that \( \eta \) is closed.
The closure of \( \eta \) together with the requirement
\( \eta|_\xi = 0 \) is equivalent to \( L_A\eta = 0 \)
for each smooth section \( A \) of~\( \xi \).

\medbreak
\nump
Given \( P \) and \( \theta \), there is a unique smooth section
\( Y \in \Gamma\zeta \) given by \( \theta(Y) = 1 \).
In particular, \( \zeta \)~is an oriented foliation.  
We have
\( L_Y\theta = Y \hook d\theta = - \pi^*\eta \).
As \( \eta \) is closed, we may locally write \( \eta = df \),
and get that
\begin{equation}
  \label{eq:loc-preservation}
  L_{e^{\pi^*f}Y}\theta = 0.
\end{equation}
Thus the vertical space \( \zeta \) for the projection
\( \pi\colon P \to M \) has local sections that preserve~\( \theta \).
This is the reason we require~\( \xi' \) to have a similar property,
since \( \xi' \) becomes the vertical spaces for the fibration
\( P \to S \) of~\( P \) over the shear~\( S \).

\medbreak
\nump
In all cases, it turns out that \( P \) carries the structure of
principal bundle.
If the fibres of \( P \) are~\( \R \),
then the flow~\( \varphi_t \) generated by~\( Y \) is a principal action.
On the other hand if we ask for the fibres of \( P \) to be circles,
the flow \( \varphi_t \) is vertical so has some
period~\( \rho(x) \) on the fibre~\( \pi^{-1}(x) \).
As \( P \) is locally trivial,
\( \rho\colon M \to \R \) is smooth,
and \( Y_0 = 2\pi Y/\pi^*\rho \) generates a vertical circle
action on~\( P \) that is principal,
cf.~\cite[Proposition 6.15]{Morita:differential-forms}.

\medbreak
\nump\label{sec:principal}
Since \( d\eta = 0 \), there is some cover \( \widetilde M \)
of~\( M \) on which (the pull-back of) \( \eta \) is exact,
\( \eta = df \).  Write \( \widetilde P \) for the pull-back
of~\( P \) to~\( \widetilde M \).

When \( P \) is an \( \R \)-bundle,
\( \tilde Y = e^fY \) generates a principal \( \R \)-action on
\( \widetilde P \) with \( \theta_0 = e^{-f}\theta \) as a
corresponding principal connection.
The curvature of \( \theta_0 \) is \( F_0 = e^{-f}F \).
As the maximal compact subgroup of \( \R \) is trivial,
\( P \) and \( \widetilde P \) are necessarily trivial and
\( [e^{-f}F] = [F_0] \) is zero in~\( H^2(\widetilde M) \).

When \( P \) is a circle bundle,
then on \( \widetilde P \) we have that the pair 
\( (\tilde Y,\tilde\theta) = (e^fY,e^{-f}\theta) \) 
satisfies \( \tilde\theta(\tilde Y) = 1 \) and
\( L_{\tilde Y}\tilde\theta = 0 \).
A priori \( \tilde Y \) is not a principal vector field as the periods
may not be constant.  
However, as we saw above there is a principal circle action which may
taken to be generated by a vector field \( Y_0 \) of the form
\( Y_0 = h\,\tilde Y \), for some positive function~\( h \)
constant on fibres.
Let \( \theta_0 \) be a principal connection with respect
to~\( Y_0 \).
Then we may write \( \tilde\theta = h\theta_0 + \nu \),
with \( \tilde Y \hook \nu = 0 \).
Now
\begin{equation}
  \label{eq:Lie-tilde}
  \begin{split}
    0
    &= hL_{\tilde Y}\tilde\theta = h\tilde Y \hook d\tilde\theta
    = Y_0 \hook d(h\theta_0 + \nu)\\
    &= - dh + Y_0\hook d\nu.
  \end{split}
\end{equation}
Choosing linearly independent one-forms~\( \alpha_i \)
on an open set \( U \subset \widetilde M \),
we may write \( \nu = \sum_{i=1}^k b_i\pi^*\alpha_i \)
with \( b_i \in C^\infty(\pi^{-1}(U)) \subset \widetilde P \).
Then \( Y_0 \hook d\nu = \sum_{i=1}^k (Y_0b_i)\pi^*\alpha_i \).
Since \( h \) is a pull-back, equation~\eqref{eq:Lie-tilde} implies 
\( Y_0b_i \) is constant on each fibre.  
As the fibres are compact, this constant must be zero, so \( dh = 0
\), and \( h \)~is constant.
We conclude that \( f \) may chosen so that
\( (\tilde Y,\tilde\theta) = (Y_0,\theta_0) \) are a principal vector
field and connection.
The curvature \( F_0 = e^{-f} F \) of \( \theta_0 = e^{-f}\theta \)
then has integral periods.

\smallbreak
\nump
As an example, consider \( \widetilde M = \R^n \setminus\{0\} \)
with \( F_0 \) a constant coefficient form with respect to the
standard coordinates \( x^i \) on~\( \R^n \).
Under the \( \R_{>0} \)-action \( m \mapsto \lambda m \),
we have \( F_0 \mapsto \lambda^2F_0 \),
so the two-form \( F = r^{-2}F_0 \),
where \( r^2 = \sum_{i=1}^n (x^i)^2 \),
is invariant.
Furthermore, \( dF = \eta \wedge F \) with
\( \eta = - d\log r^2 \).
Note that \( \log r^2 \) is not \( \R_{>0} \)-invariant.

Consider a linear isometry~\( \varphi \) of~\( \R^n \)
that also preserves~\( F_0 \).
Then for any \( \lambda > 1 \),
the group \( \Gamma \) generated by \( \lambda\varphi \)
acts freely on \( \widetilde M \) and \( M = \widetilde M/\Gamma \)
is a smooth compact manifold.
In the case that \( \varphi = \operatorname{Id} \),
\( M = S^{n-1} \times S^1 \).

On \( \R^n \), and hence \( \widetilde M \), 
\( F_0 \) is exact.
Write \( F_0 = d\alpha_0 \) with \( \alpha_0 \) homogeneous of
degree~\( 2 \) under the scaling action;
this is possible since \( F_0 \) has constant coefficients. 
We now define a principal connection~\( \theta_0 \) on the trivial
bundle~\( \widetilde P = \widetilde M \times \R \)
by \( \theta_0 = dt + \alpha_0 \).
In order that \( \theta = r^{-2}\theta_0 \) descends to a bundle
over~\( M \) we need a lift of the \( \Gamma \)-action
to~\( \widetilde P \) preserving~\( \theta \).
In the case that \( \varphi^*\alpha_0 = \alpha_0 \),
we have \( (\lambda\varphi)^*\alpha_0 = \lambda^2\alpha_0 \),
so \( (\lambda\varphi)(t) = \lambda^2t \) gives the required
lift~\( \widetilde\Gamma \).
The space \( P = \widetilde P/\widetilde G \) is now an
\( \R \)-bundle over~\( M \), with induced one-form~\( \theta \)
satisfying~\eqref{eq:d-theta} and \( \eta \) not exact.

\smallbreak
\nump
For \( P \) a circle bundle,
we claim that \( \eta \) is always exact.
Let \( \Gamma \) be the covering group of \( \widetilde M \to M \).
On \( \widetilde M \) we have \( F \) and \( \eta \) are
\( \Gamma \)-invariant.
Choose \( f \in C^\infty(\widetilde M) \) as above,
in~\S\ref{sec:principal},
with \( df = \eta \) and \( e^fY \) principal.
Invariance of \( \eta \) shows that \( g^*f-f \) is constant for each
\( g \in \Gamma \) 
and we get a homomorphism \( \lambda\colon \Gamma \to \R_{>0} \)
such that \( g^*F_0 = g^*e^{-f}F = \lambda(g)F_0 \).
The homomorphism \( \lambda \) is trivial precisely when \( \eta \) is
exact on~\( M \).
Now \( P \) is the quotient of \( \widetilde P \) by a lift of the
action of \( \Gamma \) to~\( \widetilde P \).
Furthermore, \( \theta = e^f\theta_0 \) is invariant under the lifted
action.
In particular \( g^*\widetilde P \) and \( \widetilde P \) are isomorphic
as oriented circle bundles and so the principal actions constructed
in~\S\ref{sec:principal} are isotopic.
Let \( \varphi_t \) be the isotopy from \( Y_0 \) to
\( Y_1 = g_*Y_0 \).
Then \( \varphi_1^*g^*\theta_0 \) is also a principal connection for
\( Y_0 \).
It follows that on a typical fibre \( S^1 \) of \( \widetilde P \) we
have \( 2\pi = \int_{S^1} \varphi_1^*g^*\theta_0 = \int_{S^1}
g^*\theta_0 = \lambda(g)\int_{S^1}\theta_0 = \lambda(g)2\pi \),
so \( \lambda(g) = 1 \) and \( \eta \) is exact on~\( M \).

\medbreak
\nump
It is convenient to rewrite conditions like
\eqref{eq:loc-preservation} in a global form as follows.  From
\( 0 = L_{e^gA}\theta = e^g(dg\, \theta(A) + L_A\theta)\),
under the condition \( \theta(A) \) is nowhere zero, we get 
\begin{equation}
  \label{eq:preservation}
  d\bigl(\theta(A)^{-1}L_A\theta\bigr) = 0.
\end{equation}
This is locally equivalent to the existence of a smooth
function~\( g \),
so that \( e^gA \) preserves~\( \theta \).
Another way to rewrite this is
\begin{equation}
  \label{eq:pre}
  \theta(A)\,d(A\hook d\theta)
  - d(\theta(A))\wedge (A\hook d\theta)
  = 0.
\end{equation}

Let us now apply this to the existence problem for the
foliation~\( \xi' \).
The condition \( \pi_*\xi' = \xi \),
implies that there is a unique smooth
section~\( X' \in \Gamma\xi' \)
of the form
\begin{equation}
  \label{eq:lift-X}
  X' = \widetilde X + \tilde{a}Y,
\end{equation}
where \( \widetilde X \in \Gamma\mathcal H \) is the horizontal lift
of~\( X \) to \( \mathcal H = \ker\theta \), the vector field~\( Y \)
is the unique section of~\( \zeta \) with \( \theta(Y) = 1 \) and
\( \tilde{a} \) is a smooth function.
The transversality assumption on \( \xi' \) ensures that \( \tilde{a} \)
vanishes nowhere. 
Write
\begin{equation*}
  \nu = X \hook F.
\end{equation*}
By \eqref{eq:preservation}, the condition that some local section of \( \xi' \)
preserves~\( \theta \) gives
\begin{equation*}
  \begin{split}
    0
    &= d(\tilde a^{-1}(X'\hook (\pi^*\eta\wedge\theta + \pi^*F)
    + d\tilde a))\\
    &= d(\tilde a^{-1}(-\tilde a + \pi^*\nu)) = d(\tilde a^{-1}\pi^*\nu).
  \end{split}
\end{equation*}
Contracting with \( Y \), we see that \( \tilde{a} \)~is basic 
away from zeros of~\( \nu \).
Taking \( \tilde{a} = \pi^*a \), we then have
\begin{equation}
  \label{eq:a}
  d(a^{-1} \nu) = 0.
\end{equation}
Thus a necessary condition for the existence of \( a \) is that
\begin{equation}
  \label{eq:X-F}
  d\nu \wedge \nu = 0.
\end{equation}
In other words, \( \nu = X\hook F \)~generates a differential ideal.
Note that equations \eqref{eq:a} and~\eqref{eq:X-F} are independent of
the choice of smooth generating section~\( X \) of~\( \xi \).
Using~\eqref{eq:F-eta}, equation~\eqref{eq:X-F} is seen to be
equivalent to \( L_XF\wedge (X\hook F) = 0\);
so this naturally replaces the condition in the twist
construction that \( F \) should be preserved by~\( X \).

\section{Geometric structures}
\label{sec:geometric-structures}

In the twist and shear constructions, the total space \( P \)
carries a distribution~\( \mathcal H \) that is horizontal for both
the projection~\( \pi \) to the original space~\( M \), and for the
projection to the twist or shear space.
Let us concentrate on the latter case, and denote the projection to
the shear~\( S = P/\Span{\xi'} \) by \( \pi_S\colon P\to S \).

\smallbreak
\nump
We say that two \( p \)-forms \( \alpha \in \Omega^p(M) \) and \(
\alpha_S \in \Omega^p(S) \) are \( \mathcal H \)-related,
written \( \alpha \sim_{\mathcal H} \alpha_S \), if
\( \pi^*\alpha \) and \( \pi_S^*\alpha_S \) agree on
\( \Lambda^p\mathcal H \).
This is equivalent to the existence of a \( (p-1) \)-form \(
\tilde\beta \) on~\( P \) satisfying
\begin{equation}
  \label{eq:H-rel-alpha}
  \pi^*\alpha - \pi_S^*\alpha_S = \theta \wedge \tilde\beta.
\end{equation}

\medbreak
\nump
Let us show how \( \tilde\beta \) may be specified and determine when
an \( \alpha \) gives rise to an \( \alpha_S \).
As \( \theta(Y) = 1 \), we may choose \( \tilde\beta \) so that
\( Y \hook \tilde\beta = 0 \).
Contracting \eqref{eq:H-rel-alpha} with the canonical section~\( X' \)
of~\( \xi' \) gives
\begin{equation*}
  \pi^*(X\hook \alpha)
  = \pi^*a\, \tilde\beta - \theta\wedge(\widetilde X \hook \tilde\beta). 
\end{equation*}
The horizontal part of this equation says
\( \tilde\beta = \pi^*\beta \) for \( \beta = a^{-1}X\hook\alpha \).
This gives \( \widetilde X \hook \pi^*\beta = 0 \), and so the
vertical part of the equation is also satisfied.
Differentiating~\eqref{eq:H-rel-alpha}, we get
\begin{equation}
  \label{eq:d-alpha-hor}
  \pi^*d\alpha - \pi_S^*d\alpha_S
  = \pi^*(F \wedge \beta) - \theta \wedge \pi^*(\eta\wedge\beta + d\beta).
\end{equation}
Again contracting with \( X' \), we get horizontally
\begin{equation*}
  \begin{split}
    X \hook d\alpha
    &= \nu\wedge\beta - a(d\beta + \eta\wedge\beta) \\
    &= \bigl(a^{-1}(da + \nu) - \eta\bigr)\wedge (X\hook \alpha)
    - d(X\hook \alpha),
  \end{split}
\end{equation*}
which implies the vanishing of the vertical component too.
In other words, \( \alpha \)~satisfies
the \emph{automorphic} condition
\begin{equation}
  \label{eq:invariance}
  L_X\alpha = \gamma \wedge (X\hook \alpha),\quad
  \text{for}\ 
  \gamma = a^{-1}(da + \nu) - \eta.
\end{equation}
Note that invariant functions and \( \eta \) are automorphic.
Also \eqref{eq:a} implies that \( F \) is automorphic in
this sense.
The defining property of~\( a \) ensures that \( \gamma \) is closed,
so locally we may write \( \gamma = dh \) and \eqref{eq:invariance} becomes
\begin{equation}
  \label{eq:invar-h}
  L_{e^{-h}X}\alpha = 0,\quad \text{for}\ dh = \gamma.
\end{equation}

For \( \alpha \) automorphic, \( \mathcal H \)-related
to~\( \alpha_S \), equation~\eqref{eq:d-alpha-hor} shows that
\begin{equation}
  \label{eq:dS}
  d\alpha_S \sim_{\mathcal H} d\alpha - a^{-1}F \wedge (X\hook \alpha)
  \eqqcolon d_S\alpha.
\end{equation}
It follows that the automorphic forms on~\( M \) form a
differential algebra under~\( d_S \).

\smallbreak
\nump
In the twist construction \( \eta = 0 \) and the only forms considered
for transferring to the twist space where invariant forms.
Thus in condition~\eqref{eq:invariance} \( \gamma \) is zero,
which translates to \( da = - X\hook F \), agreeing with the
definition of~\( a \) given in \cite{Swann:twist}.
In this case formula \eqref{eq:dS} agrees with that in the twist
construction. 

In the algebraic model of the shear, the function \( a \) is constant, 
equal to~\( -1 \), and all left-invariant forms are automorphic

\section{Examples}
\label{sec:examples}

\nump
Consider the Abelian algebra \( \lie g = \R^6 \) with an invariant
K\"ahler structure \( (g,\omega,I) \).
We may change bases so that \( \omega = e_{12} + e_{34} + e_{56} \),
the \( (1,0) \)-forms for~\( I \) are
\( \Lambda^{1,0} = \Span{e_1 + ie_2,e_3 + ie_4, e_5 + ie_6} \)
and a given symmetry is \( X = E_1 \).
Let us consider a shear for which each \( e_i \) is automorphic.
As \( L_Xe_i = 0 \), this means
\( 0 = \gamma = a^{-1}(da + X\hook F) - \eta \).
With this satisfied, \( g = \sum_{i=1}^6 e_i^2 \) and \( \omega \) are
automorphic.

Let us find shears of \( \lie g \) for which
the induced geometry is K\"ahler.
The closure of \( \omega_S \) gives
\( 0 = d_S\omega = -a^{-1}F\wedge e_2 \),
so \( F = e \wedge e_2 \) for some one-form~\( e \).
Integrability of the shear~\( I_S \) of~\( I \)
is equivalent to
\( d_S\Lambda^{1,0} \subset \Lambda^{2,0} + \Lambda^{1,1} \).
However, \( d_Se_i = 0 \) for \( i>1 \), so the only condition comes
from \( d_S(e_1 + ie_2) = - a^{-1}F \) and implies that \( F \) is of
type \( (1,1) \).
We conclude that \( F = be_{12} \) for some function~\( b \).
For this example, we take \( b = 1 \).
Now the condition \( dF = \eta \wedge F \) with \( d\eta = 0 \) and
\(\eta(X) = 0 \) gives \( \eta = \mu e_2 \) for some smooth
function~\( \mu(y) \) of the second coordinate~\( y \) on~\( \R^6 \).
Again for simplicity we take \( \mu = 1 \).
Now \( d(a^{-1}X\hook F) = 0 \), implies \( a = a(y) \) and
\( \gamma = 0 \) reduces to \( a_y = a - 1 \), 
which has solutions \( a = ce^y + 1 \).
For \( c = 0 \), the resulting K\"ahler shear is the solvable algebra
\( \lie r = (12, 0^5) = (12,0) \oplus \R^4 \).
For \( c\ne 0 \), we get other examples which are defined
for all \( y \) if \( c > 0 \); smoothly these are still the space
\( \mathcal H^2\times \R^4 \), where \( \mathcal H^2 \) is the
hyperbolic plane, but now with a non-invariant K\"ahler structure on the
first factor.

\medbreak
\nump
Consider the almost Abelian 
\begin{equation*}
 \lie g_{\lambda,\mu} = ((\lambda+\mu) 17, \lambda 27, \mu 37,
-(\lambda+\mu) 47, -\lambda 57,-\mu 67,0)
\end{equation*}
where \( \lambda \) and \( \mu \) are real constants.
By \cite{Freibert:cocal} this carries a co-calibrated
\( G_2 \)-structure with closed four-form
\( \psi = e_{1425} + e_{1436} + e_{2536} - e_{4567} + e_{4237} +
e_{1267} + e_{1537} \).

Let \( X = E_1 \) and \( F = e_{23} \).
Then \( \nu = X \hook F = 0 \), so the condition~\eqref{eq:a}
is satisfied for any function~\( a \).
We have \( dF = \eta\wedge F \) with \( \eta = -(\lambda+\mu)e_7 \).
For \( i>1 \), each \( e_i \) has both \( L_Xe_i = 0 \) and
\( X\hook e_i = 0 \), so is automorphic.
The automorphic condition for \( e_1 \), is that 
\( L_Xe_1 = X \hook de_1 = (\lambda+\mu)e_7 = - \eta \)
should equal
\( \gamma \wedge X\hook e_1 = \gamma = a^{-1}da - \eta \),
which forces \( a \) to be constant.
It follows from \eqref{eq:dS} that the shear of \( \lie g_{\lambda,\mu} \) is an
algebra \( \lie h_{\lambda,\mu} \) whose dual is generated by \(
(e_1)_S,\dots,(e_7)_S \), the one-forms \( \mathcal H \)-related to \(
e_1,\dots,e_7 \) and with differentials given by
\( ((\lambda+\mu) 17 - a^{-1} 23, \lambda 27, \mu 37,
-(\lambda+\mu) 47, -\lambda 57,-\mu 67,0)) \).
These algebras \( \lie h_{\lambda,\mu} \) are no longer almost Abelian.
The form \( \psi_S \) that is \( \mathcal H \)-related to~\( \psi \)
is also closed, since \( d\psi_S \) is \( \mathcal H \)-related to
\( d_S\psi = d\psi - a^{-1}F\wedge X\hook \psi
= 0 - a^{-1}e_{23}\wedge(e_{425}+e_{436}+e_{267}) = 0 \).
Thus the shear induces a co-calibrated \( G_2 \)-structure on
each~\( \lie h_{\lambda,\mu} \). 

\medbreak
\nump
Let us more generally consider shears that preserve geometries with
closed forms.
Equation~\eqref{eq:dS} implies that if \( \alpha \) is a closed
\( p \)-form, then the form \( \alpha_S \) \( \mathcal H \)-related to
\( \alpha \) is closed if and only if \( F\wedge X\hook \alpha = 0 \).
Thus for \( \omega \) a symplectic form, \( \omega_S \) closed implies
that \( F \) is necessarily decomposable, much as we saw in the
previous example.

In dimension~\( 6 \), suppose we have an \( \SU3 \)-structure
\( (g,J,\omega,\rho_+,\rho_-) \),
where \( g \) is a Riemannian metric,
\( J \)~a compatible almost complex structure,
\( \omega(\cdot,\cdot) = g(J\cdot,\cdot) \) and
\( \rho_+ + i\rho_- \) is a compatible complex volume form.

This is half-flat if (i)~\( \omega \) is co-symplectic, meaning
\( d(\omega^2) = 0 \), and (ii)~\( d\rho_- = 0 \).
The shear form \( \omega_S \) is co-symplectic if away from zeros of
\( X \) we have
\( 0 = F\wedge X\hook \omega^2 = 2F\wedge JX^\flat \wedge \omega \).
This holds if and only if \( F = \beta \wedge JX^\flat + F_0 \),
where \( \beta \) is any one-form and \( F_0 \in \Lambda^2W \),
\( W = \Ann\Span{X,JX} \) the annihilator of \( X \) and~\( JX \),
is a two-form orthogonal to \( \omega_0 = \omega - g(X,X)^{-1}X^\flat\wedge
JX^\flat \in \Lambda^2W \).
For the shear form \( (\rho_-)_S \) to be closed, we need
\( 0 = F\wedge X\hook \rho_- \).
Now \( \rho_- \) is proportional to
\( X^\flat\wedge \omega_1 + JX^\flat\wedge\omega_2 \)
with \( \omega_1 + i\omega_2 \) a complex symplectic form in
\( \Lambda^{2,0}W \).
Indeed \( \omega_0,\omega_1,\omega_2 \) span the space \( \Lambda^2_-W
\) of anti-self dual two-forms in \( \Lambda^2W \).
Thus we require \( F \wedge \omega_1 = 0 \).
Together with the co-symplectic condition we find that
\( F \in \Lambda^2_+W + \R\omega_2 \).

Let \( (N^4,\omega_I,\omega_J,\omega_K) \) be a K3 surface endowed
with a hyperK\"ahler structure.  
Taking \( M = N^4 \times S^1\times S^1 \), we obtain an \( \SU3
\)-structure by putting \( \omega = \omega_I + dx\wedge dy \),
\( \rho_+ + i\rho_- = (\omega_J+i\omega_K)\wedge(dx + idy) \).
Let \( X = \partial_x \), and consider the cohomology of~\( N \).
As is well-known, see~\cite{Besse:Einstein}, \( H^2(N) = H^2_+(N)
\oplus \Span{[\omega_I],[\omega_J],[\omega_K]} \) with
\( b^2_+(N) = 19 \).
Let \( G \in \Omega^2_+(N) \) be a harmonic representative for an
primitive integral element of~\( H^2_+(N) \),
and consider \( F = e^fG \) for \( f \in C^\infty(M) \).
Now \( dF = \eta\wedge F \), with \( \eta = df \).
The condition \( \eta|_\xi = 0 \), implies \( \partial_xf = 0 \).
We have that \( F \) is orthogonal to \( \omega_I \), \( \omega_J \)
and \( \omega_K \),
so shearing with \( (X,F,\eta) \) the \( \SU3 \)-structure remains 
half-flat.
Note that as \( \nu = X\hook F = 0 \),
we may choose \( a \) so that \( \omega \) and \( \rho_- \)
are automorphic, since both forms are \( X \)-invariant and taking
\( a = Ce^f \), \( C \)~constant, 
gives \( \gamma = 0 \) in~\eqref{eq:invariance}.
We thus obtain half-flat \( \SU3 \)-structures on spaces
\( Q^5 \times S^1 \), where \( Q \to N \) is a circle-bundle 
with curvature~\( -C^{-1} G \).
The topology constrains~\( 1/C \) to be an integer.

In this example, the automorphic condition for \( \rho_- \) is
\( L_X\rho_- = 0 = \gamma \wedge X\hook\rho_- = \gamma \wedge
g(X,X)\omega_J \).  As \( \omega_J \) has rank~\( 4 \), this implies
\( \gamma = 0 \) and \( a^{-1}da = \eta \). 
On the other hand the automorphic condition for \( \omega^2 \) is
weaker:
\( L_X\omega^2 = 0 = \gamma \wedge X\hook \omega^2 = 2\gamma\wedge
JX^\flat \wedge \omega = 2\gamma \wedge dy \wedge \omega_I \),
this is equivalent to \( \gamma \wedge dy = 0 \).
However \( \gamma \) is closed, so \( \gamma = dH \) with \( H=H(y) \)
just a (periodic) function of~\( y \). 
So if looking for just co-symplectic shears, we thus have
\( a^{-1}da = \eta + dH \), so \( a \)~is the more general function
\( a = Ce^{f+H(y)} \).

\medbreak
\nump
Let \( C = \R \times Q \) with \( Q \) a principal circle bundle over
a base~\( B \).
Write \( \nu \) for the principal connection and
\( \omega = d\nu \) for its curvature.
Consider \( F = dt \wedge \nu + \omega \),
where \( t \) is the standard coordinate on~\( \R \).
We have \( dF = -dt \wedge \omega = \eta \wedge F \),
for \( \eta = -dt \).
Let \( X \) generate the principal action on~\( Q \),
then \( X\hook F = -dt \) and equation~\eqref{eq:a} is satisfied by
any function~\( a(t) \).
The one-form \( \gamma \) in the automorphic
condition~\eqref{eq:invariance} is proportional to~\( dt \).
In particular, any invariant form~\( \alpha \) with \( X\hook\alpha \)
in the ideal generated by \( dt \) is automorphic.
But in general it is important to consider non-invariant automorphic
forms as the Lie algebra examples emphasise.

\providecommand{\bysame}{\leavevmode\hbox to3em{\hrulefill}\thinspace}

\end{document}